%% file: main.tex
\documentclass[letterpaper, 10 pt, conference]{ieeeconf}  % Comment this line out
\IEEEoverridecommandlockouts                              % This command is only
\input{packages.tex}

\input{commands.tex}

\title{\LARGE \bf
Opinion Dynamics Steering using Stochastic Search
}

\author{Ziyi Wang% <-this % stops a space
\thanks{Ziyi Wang is with Center for Machine Learning,
        Georgia Institute of Technology, Atlanta, GA, USA.} 
\thanks{Evangelos A. Theodorou is with Faculty of Aerospace Engineering Department, Georgia Institute of Technology, Atlanta, GA, USA.}\thanks{Correspondence to: \href{mailto: ZiyiWang@gatech.edu}{ZiyiWang@gatech.edu}} and Evangelos A. Theodorou}

\begin{document}

\maketitle
\thispagestyle{empty}
\pagestyle{empty}

%%%%%%%%%%%%%%%%%%%%%%%%%%%%%%%%%%%%%%%%%%%%%%%%%%%%%%%%%%%%%%%%%%%%%%%%%%%%%%%%
\begin{abstract}

In this paper, we apply the stochastic search dynamic optimization framework for steering opinion dynamics in a partially active population. The framework is grounded on stochastic optimization theory and relies on sampling of candidate solutions from distributions of the exponential family. We derive the distribution parameter update for an open loop, a feedback, and a novel adaptive feedback policy. All three policies are tested on two different opinion dynamics steering scenarios in simulation and compared against a hand designed feedback policy. The results showcase the effectiveness the framework for opinion dynamics control and the robustness of the adaptive feedback policy with respect to the active agent set size.

\end{abstract}

%%%%%%%%%%%%%%%%%%%%%%%%%%%%%%%%%%%%%%%%%%%%%%%%%%%%%%%%%%%%%%%%%%%%%%%%%%%%%%%%
\section{Introduction}
\label{sec: intro}

There has been a wide spectrum of interdisciplinary research efforts aimed to understand and model the complex human behavior for decades. Sociologists and psychologists study the formulation of opinions and social preferences \cite{little1952social, french1956formal}. Economists investigate the cascade of information and how individual decisions lead to complex economic phenomena \cite{bikhchandani1992theory, puu2013attractors}. In recent years, there is also a growing interest from politicians to measure and predict the public opinion \cite{nall2015political, lamberson2018model}. The ultimate goal of relevant research in different fields is to influence and control the underlying human opinion process.

There exists extensive prior work on the modeling of opinion dynamics, dating back to the seminal works of DeGroot and Friedkin \cite{degroot1974reaching, friedkin1990social}. Since then, a variety of both theoretical and application-oriented models have been proposed \cite{olfati2004consensus, fortunato2005vector, quattrociocchi2014opinion, friedkin2011social, xu2020paradox}. However, despite the abundance of literature on the modeling side, research on opinion dynamics control remains limited and mainly focuses on pinning control for deterministic opinion dynamics models \cite{song2010second, kan2015containment, delellis2017steering}. Pinning control establishes a virtual leader that can influence the agents in a network. The framework designs a feedback policy on the error vector between the leader and the followers to track the desired trajectory of the leader.

In the meantime, stochastic optimal control and stochastic optimization algorithms have achieved much success in the field of robotics and autonomy \cite{williams2017information, williams2017model}. In this work, we tackle the problem of opinion dynamics steering from the perspective of stochastic optimization theory. We leverage the stochastic search framework \cite{zhou2014gradient} that solves general nonlinear stochastic optimization problems. The framework relies on relaxation of the original problem and samples candidate solutions from parameterized distributions of the exponential family. The distribution parameter updates are derived via gradient descent. Many well-known algorithms, such as \ac{CEM} \cite{de2005tutorial}, Evolutionary Strategies (ES) \cite{hansen2016cma}, and \ac{MPPI} \cite{williams2017information}, are special cases of the framework. Due to its generality and ability to handle arbitrarily complex dynamics and cost function, stochastic search has been applied to many dynamic optimization problems in the last few years, including constrained trajectory optimization \cite{boutselis2020constrained}, jump diffusion process control \cite{wang2019information}, risk sensitive control \cite{wang2021adaptive}, and even multi-agent swarm control \cite{fan2018model}.

In this paper, we apply the stochastic search framework to the opinion dynamics model proposed in \cite{xu2020paradox}, which is a stochastic variant of the well-known Friedkin and
Johnsen model \cite{friedkin1990social}. We consider the problem of opinion dynamics steering for partially actuated population in which only a subset of all agents in the population can be controlled. We derive the update law for an open loop and a feedback policy. In addition, we design an adaptive feedback policy that optimizes with respect to which agents to actuate while keeping the active agent set size constant in expectation. The proposed framework is tested in simulation on 2 different opinion steering tasks for a network of 200 agents. 

The main contribution of this paper is threefold:
\begin{enumerate}
    \item We derive the stochastic search policy update for the opinion dynamics model.
    \item We design a novel adaptive feedback policy that simultaneously optimizes the feedback parameters and the active agent set.
    \item We showcase the opinion dynamics steering capability of the algorithm in 2 different scenarios and compare against a hand designed baseline policy. The results also demonstrate the robustness of the adaptive feedback policy with respect to the active agent set size.
\end{enumerate}

The rest of this paper is organized as follows: sections \ref{sec: notation} through \ref{sec: stochastic_search} introduce the notation, opinion dynamics model, and stochastic search dynamic optimizer. In \cref{sec: ss_control}, we derive the stochastic search policy update. We present the adaptive stochastic search algorithm in \cref{sec: algo}. The simulation results are included in \cref{sec: simulation}. Finally, we conclude the paper and discuss future directions in \cref{sec: conclusion}.

%%%%%%%%%%%%%%%%%%%%%%%%%%%%%%%%%%%%%%%%%%%%%%%%%%%%%%%%%%%%%%%%%%%%%%%%%%%%%%%%
\section{Notation}
\label{sec: notation}
Throughout the paper, we use subscript/superscript $i$ to denote the state of the $i$th agent in a population. Without subscript/superscript the variable refers to the collection of the state across all agents $x_t=\{x_t^1, x_t^2, \cdots, x_t^N\}$. Upper case letters are used to represent trajectories of the variable over the time horizon $X=\{x_1, x_2, \cdots, x_T\}$. For notational simplicity, the distributions that expectations are taken with respect to are only written explicitly the first time they appear.

\section{Opinion Dynamics Model}
\label{sec: opinion_model}

\begin{figure}[t!]
    \centering
    \includegraphics[width=0.9\linewidth]{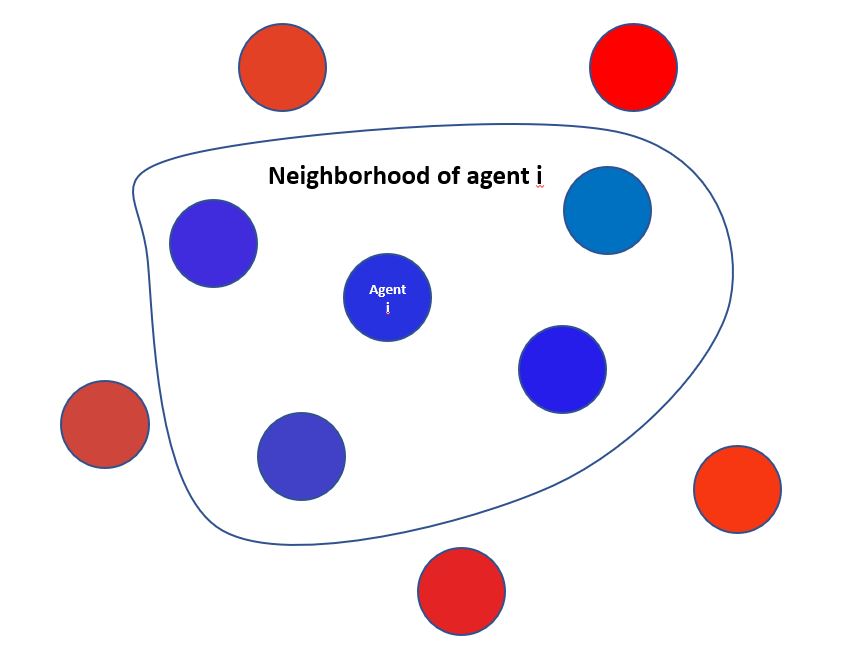}
    \caption{Schematic of neighborhood of agent i. The color indicates the state of each agent and the neighborhood consists of agents whose state is similar to the ego agent.}
    \label{fig: neighborhood}
\end{figure}
In this section we introduce the opinion dynamics model \cite{xu2020paradox}. Consider a group of $N$ agents with each agent having the state (opinion) $x^i_t$ over a discrete time horizon $T$. The dynamics of each agent is modeled as follows:
\begin{equation} \label{eq: initial_dynamics}
x^i_{t+1} = (1-\alpha) x^i_t + \alpha f(\calX^i_t) + \sigma w^i_t,
\end{equation}
where the first two terms denote the drift and the last term represent the state fluctuations in the form of independent Gaussian noise. Here $\alpha\in[0,1]$ is the \textit{susceptibility to local influence}. The function $f(\calX^{i}(t))$ is the \textit{center of bias} and characterizes how the neighborhood affects the agent's opinion. The neighborhood $\calX^i_t$ and bias function $f$ are defined as:
\begin{equation} \label{eq: neighborhood}
\calX_t^i = \{x^j_t, \forall j\in\{1,\cdots,N\}\setminus i\, |\, \|x^i_t - x^j_t\|_2 \leq \epsilon\},
\end{equation}
and
\begin{equation} \label{eq: bias_function}
f(\calX_t^i) = \frac{1}{|\calX^i_t|}\sum_{j\in\calX_t^i} x_j.
\end{equation}
Here $|\calX^i_t|$ denotes the size of the neighborhood set. The neighborhood, as illustrated in \cref{fig: neighborhood}, is defined as a 2-norm ball around the agent's current state that "pulls" the agent towards its mean state.

In this work we consider an actuation model that directly influences the agent's state in the form of $u^i_t$ added to \eqref{eq: initial_dynamics}. In most practical scenarios, however, not all agents in a population can be actuated. In this case we model two sets of agents within the same population, namely a set of passive agents $ I_{P} $  the opinion of which is only influenced only by other agents and the active set $I_{A} $ of agents whose state can be actuated based on an objective function. In mathematical terms this leads to the following representation:
\begin{align}\label{eq: passive_dynamics}
x^i_{t+1} &= (1-\alpha) x^i_t + \alpha f(\calX^i_t) + \sigma w^i_t, \quad\quad \, \ \forall i\in I_P\\
x^i_{t+1} &= (1-\alpha) x^i_t + \alpha f(\calX^i_t) + u^i_t + \sigma w^i_t, \ \forall i\in I_A.
\label{eq: active_dynamics}
\end{align}

% \addtolength{\textheight}{-3cm}   % This command serves to balance the column lengths
%                                   % on the last page of the document manually. It shortens
%                                   % the textheight of the last page by a suitable amount.
%                                   % This command does not take effect until the next page
%                                   % so it should come on the page before the last. Make
%                                   % sure that you do not shorten the textheight too much.

% %%%%%%%%%%%%%%%%%%%%%%%%%%%%%%%%%%%%%%%%%%%%%%%%%%%%%%%%%%%%%%%%%%%%%%%%%%%%%%%%
\section{Stochastic Search}
\label{sec: stochastic_search}
In this section we introduce the stochastic search framework for solving general dynamic optimization problems of the form:
\begin{equation}\label{eq: ss_initial_obj}
\begin{split}
\pi^* &= \argmin_\pi \Eb_{p(\tau)}[J(\tau)]\\
\text{s.t.}\hspace{3mm} \ x_{t+1} &= F(x_t, \pi(x_t;\phi_t), w_t)\\
x_0 & \sim p(x_0),
\end{split}
\end{equation}
where $\tau=\{X,U\}$ denote the state-control trajectory pair and $w_t$ is the stochasticity in dynamics. Note that both the cost function and dynamics can be arbitrarily complex and even discontinuous. The controls are generated from a policy $u_t=\pi(x_t;\phi_t)$ parameterized by $\phi_t$. Examples of such policies are 1) open loop policy, $u_t = \phi_t$; 2) feedback policy, $u_t = K_t x_t + k_t, \phi_t = \{K_t, k_t\}$; 3) neural network policy with weights $\phi_t$.

\subsection{Policy Update Derivation}

We now assume that the policy parameters are sampled from a distribution of the exponential family $\phi_t\sim p(\phi_t;\theta_t)$ such that 
\begin{equation} \label{eq: exp_family_dist}
p(\phi_t;\theta_t) = h(\phi_t)\exp\left(\eta(\theta_t)\T T(\phi_t) - A(\theta_t)\right),
\end{equation}
where $T(\phi_t)$ denotes the sufficient statistic of the distribution, and $\eta(\theta_t)$ denote the natural parameters, which are functions of the distribution parameters $\theta_t$. The optimization problem is now transformed into
\begin{equation} \label{eq: ss_exp_obj}
\Theta^* = \argmin_\Theta \Eb_{p(\tau;\Theta)}\left[J(\tau; \Theta)\right],
\end{equation}
where $\Phi$ and $\Theta$ are the policy and distribution parameter trajectory. Note that the new objective is an upper bound on the original objective $\Eb_{p(\tau;\Theta)}\left[J(\tau; \Theta)\right]\geq \Eb_{p(\tau)}[J(\tau)]$ and equality is achieved when all probability mass is on deterministic solution.

To be consistent with literature on stochastic search in optimization, we turn the minimization problem into a maximization one by optimizing with respect to $-J$. We then apply an non-decreasing shape function $S:\Rb\rightarrow \Rb^+$ to the cost function. Due to its monotonicity, applying the shape function does not change the solution set of the problem but provide added flexibility in algorithm design. Common shape functions include: 1) the exponential function, $S(y;\lambda)=\exp(\lambda y)$, which leads to the \ac{MPPI} update law \cite{williams2017information}; 2) the sigmoid function, $S(y; \lambda, \psi)=(y-y_{\text{lb}})\tfrac{1}{1+\exp(-\lambda(y-\psi))}$, where $y_{\text{lb}}$ is a lower bound for the cost and $\psi$ is the ($1-\psi$)-quantile, which results in an update law similar to the \ac{CEM} \cite{de2005tutorial} but with soft elite threshold $\psi$. The sharpness of the cutoff is controlled by $\lambda$.

Finally, we apply another log transformation to obtain a gradient invariant to the scale of the
objective function. The optimization problem becomes
\begin{equation}
\Theta^* = \argmax_\Theta \log \Eb\left[S(-J(\tau;\Theta))\right] = l(\Theta).
\end{equation}
To derive the parameter update, we first decompose the path probability of trajectories as
\begin{equation}
p(\tau;\Theta) = p(x_0)\prod_{t=0}^{T-1} p(x_{t+1}|x_t, \pi(x_t;\phi_t)) p(\phi_t;\theta_t),
\end{equation}
where $p(x_{t+1}|x_t, \pi(x_t;\phi_t))$ is the probabilistic representation of the stochastic dynamics. With this we can compute the gradient of the objective with respect to the natural parameter at each timestep as
\begin{align}
\nabla_{\eta_t} l(\Theta) &= \frac{\nabla_{\eta_t} \Eb\left[S(-J(\tau;\Theta))\right]}{\Eb\left[S(-J(\tau;\Theta))\right]}\\
&= \frac{\Eb\left[S(-J(\tau;\Theta))\nabla_{\eta_t}\ln p(\phi_t;\eta_t)\right]}{\Eb\left[S(-J(\tau;\Theta))\right]}\\
&= \frac{\Eb\left[S(-J(\tau;\Theta))(T(\phi_t) - \Eb[T(\phi_t)])\right]}{\Eb\left[S(-J(\tau;\Theta))\right]},
\end{align}
where the "log trick", $\nabla p=p\nabla\ln p$, and the Markovian property of the dynamics are used. The natural parameters can be updated iteratively through gradient ascent:
\begin{equation} \label{eq: ss_update}
\eta^{k+1}_t = \eta^k_t + \beta^k \frac{\Eb\left[S(-J(\tau;\Theta))(T(\phi_t) - \Eb[T(\phi_t)])\right]}{\Eb\left[S(-J(\tau;\Theta))\right]}.
\end{equation}
The distribution parameter update can then be computed by plugging in the definition of $\eta_t$ as a function of $\theta_t$.

\section{Stochastic Search for Opinion Dynamics Control}
\label{sec: ss_control}
In this section, we formulate the opinion steering problem and apply the stochastic search framework for 3 different policies. The opinion steering problem can be formulated as:
\begin{equation} \label{eq: opinion_steering_problem}
\begin{split}
& \ \pi^* = \argmin_\pi \Eb\left[\sum_{t=1}^{T}\left(\sum_{i=1}^N q\|x^i_t - \zeta\|_2^2 + \sum_{i\in I_A}r\|u^i_t\|^2_2\right)\right]\\
\text{s.t.} & \ x^i_{t+1} = (1-\alpha) x^i_t + \alpha f(\calX^i_t) + \sigma w^i_t, \quad\quad \, \ \forall i\in I_P\\
& \ x^i_{t+1} = (1-\alpha) x^i_t + \alpha f(\calX^i_t) + \pi(x^i_t;\phi_t) + \sigma w^i_t, \ \forall i\in I_A\\
& \ x^i_0 \sim p(x_0).
\end{split}
\end{equation}
The goal is to steer the partially controlled population from an initial distribution to a target state $\zeta$. Since the problem formulation fits the form considered in \eqref{eq: ss_initial_obj}, the stochastic search framework can be readily applied to the problem with iterative update law derived in \eqref{eq: ss_update}.
\subsection{Open loop and feedback Policy}
\label{sec: open_feedback_policy}
With the generic update law for distribution's natural parameters, we can derive specific parameter updates for an open loop, $u^i_t = \phi^i_t, \phi_t = \{\phi^i_t\}_{i=I_A}$, and a feedback policy, $u^i_t = K_x x^i_t + k_t, \phi_t = \{K_t, k_t\}$. Note that the same feedforward and feedback terms are applied to all controlled agents. In both cases, we assume that the policy parameters are sampled from a Gaussian distribution $\calN(\mu,\Sigma)$ with fixed variance and varying mean. In this case the distribution parameter is only the mean $\mu$. The natural parameters and sufficient statistic of the distribution are $\eta(\theta)=\mu\Sigma^{-\frac{1}{2}}$ and $T(\phi)=\phi\Sigma^{-\frac{1}{2}}$ respectively. Plugging these into \eqref{eq: ss_update} we get:
\begin{equation} \label{eq: gaussian_update}
\begin{split}
\mu_t^{k+1}\Sigma_t^{-\frac{1}{2}} &= \mu^k_t\Sigma_t^{-\frac{1}{2}} + \beta^k \frac{\Eb\left[S(-J)(\phi_t\Sigma^{-\frac{1}{2}} - \Eb[\phi_t\Sigma^{-\frac{1}{2}}])\right]}{\Eb\left[S(-J)\right]}\\
\Rightarrow \mu^{k+1}_t &= \mu^k_t + \beta^k \frac{\Eb\left[S(-J)(\phi_t - \mu^k_t)\right]}{\Eb\left[S(-J)\right]}.
\end{split}
\end{equation}
The expectation can be approximated via Monte Carlo sampling as:
\begin{equation} \label{eq: gaussian_numerical_update}
\mu^{k+1}_t = \mu^k_t + \beta^k \frac{\sum_{m=1}^M S(-J^m)(\phi^m_t - \mu^k_t)}{\sum_{m'=1}^M S(-J^{m'})}.
\end{equation}

\subsection{Adaptive feedback Policy}
\label{sec: adaptive_policy}
Alternatively, we can formulate an adaptive feedback policy that also optimizes with respect to which agents are actuated. Assume that the percentage of actuated agents $\xi$ is fixed $|I_A| = \xi N$. The active agent set $I_A$ is now part of the optimization. We assign an actuation indicator variable $a_i$ to each agent in the population that represents whether the agent is actuated. Each variable is sampled from the Bernoulli distribution and has an associated parameter $p_i$ that denotes the probability of actuation. Since the active agent set changes throughout optimization, a feedback policy $u^i_t=K_t x^i_t$ is applied to the actuated agents as opposed to the open loop policy.

The optimization is then performed with respect to both the policy parameters $\phi$ and actuation probability $p$ as
\begin{equation} \label{eq: adaptive_opinion_steering_problem}
\begin{split}
\pi^*, p^* = &\argmin_{\pi, p} \Eb\left[\sum_{t=1}^{T}\left(\sum_{i=1}^N q\|x^i_t - \zeta\|_2^2 + \sum_{i\in I_A}r\|u^i_t\|^2_2\right)\right]\\
\text{s.t.} & \ x^i_{t+1} = F(x^i_t, \pi(x^i_t;\phi_t), p^i, w^i_t)\\
& \ x^i_0 \sim p(x_0).
\end{split}
\end{equation}
The feedback policy parameters are again assumed to be sampled from a Gaussian distribution and can be updated using \eqref{eq: gaussian_numerical_update}. For the Bernoulli random variables representing the actuation probability, the natural parameters and sufficient statistic of the distribution are $\eta(\theta)=\log\frac{p}{1-p}$ and $T(a) = a$. Plugging these into the parameter update \eqref{eq: ss_update} we get:
\begin{equation}
\begin{split}
\log &\frac{p^{k+1}_i}{1-p_i^{k+1}} = \log \frac{p^{k}_i}{1-p_i^{k}} + \beta^k \frac{\Eb[S(-J)(a_i-\Eb[a_i])]}{\Eb[S(-J)]}\\
\Rightarrow \log& \frac{1-p_i^{k+1}}{p_i^{k+1}} = \log \frac{1-p_i^{k}}{p_i^{k}} - \beta^k \frac{\Eb[S(-J)(a_i-p_i^k)]}{\Eb[S(-J)]}\\
\Rightarrow &\frac{1}{p_i^{k+1}} - 1 = (\frac{1}{p_i^k} - 1)\exp\left(- \beta^k \frac{\Eb[S(-J)(a_i-p_i^k)]}{\Eb[S(-J)]}\right)
\end{split}
\end{equation}
which leads to the actuation probability update:
\begin{equation} \label{eq: bernoulli_update}
p^{k+1}_i =
\left(1 + (\frac{1}{p_i^k} - 1)\exp\left(- \beta^k \frac{\Eb[S(-J)(a_i-p_i^k)]}{\Eb[S(-J)]}\right)\right)^{-1}
\end{equation}
The expectation can also be approximated via sampling as:
\begin{equation} \label{eq: bernoulli_numerical_update}
\begin{split}
&p^{k+1}_i =\\
&\left(1 + (\frac{1}{p_i^k} - 1)\exp\left(- \beta^k \frac{\sum_{m=1}^M S(-J^m)(a^m_i-p_i^k)}{\sum_{m'=1}^M S(-J^{m'})}\right)\right)^{-1}
\end{split}
\end{equation}
Since the actuation probability update \eqref{eq: bernoulli_numerical_update} and policy parameter update \eqref{eq: gaussian_numerical_update} are independent, they can be performed simultaneously during each optimization iteration.

\begin{algorithm}[t]
\caption{Adaptive Feedback Stochastic Search Controller}
\begin{algorithmic}
    \STATE \textbf{Given:} susceptibility: $\alpha$;  noise level: $\sigma$; neighborhood size: $\epsilon$; active portion: $\xi$; planning horizon $T$; MPC horizon: $T_{\text{mpc}}$; initial distribution: $p(x_0)$; Number of samples: $M$; Number of agents: $N$; Sampling variance $\Sigma$; Step size: $\beta$;
    \STATE Sample initial state $x_0\sim p(x_0)$;
    \STATE Initialize policy $p^0, \mu^0$;
    \FOR{$t_{\text{mpc}}=0$ to $T_{\text{mpc}}$}
    \item Copy initial state $\{\hat{x}^i_0\}_{i=1,\cdots,N}=x^i_{t_{\text{mpc}}}$;
    \FOR{$m=1$ to $M$}
    \FOR{$i=1$ to $N$}
    \item Sample $a^{i.m}\sim\text{Bernoulli}(p^{i,t_{\text{mpc}}})$;
    \FOR{$t=0$ to $T-1$}
    \item Sample $K^m_t\sim \calN(\mu_t^{t_{\text{mpc}}}, \Sigma)$;
    \item Sample $w^m_t\sim\calN(0, \sigma^2)$
    \item Compute control $u_t^{i,m}=K_t^m \hat{x}^i_t$;
    \IF{$a^{i,m}=0$}
    \item Propagate passive dynamics $\hat{x}^i_{t+1}$ \eqref{eq: passive_dynamics};
    \ELSE
    \item Propagate active dynamics $\hat{x}^i_{t+1}$ \eqref{eq: active_dynamics};
    \ENDIF
    \ENDFOR
    \ENDFOR
    \item Compute cost $J^m(\hat{X}^m, U^m)$;
    \ENDFOR
    \item Update $\mu^{t_{\text{mpc}}+1}$ \eqref{eq: gaussian_numerical_update};
    \item Update $\bar{p}^{t_{\text{mpc}}+1}$ \eqref{eq: bernoulli_numerical_update};
    \item Normalize actuation probability $p^{t_{\text{mpc}}+1}=\frac{\xi N\bar {p}^{t_{\text{mpc}}+1}}{\sum_{i=1}^N \bar {p}^{t_{\text{mpc}}+1}_i}$;
    \item Execute control $\mu^{t_{\text{mpc}}+1}_0$;
    \FOR{$t=0$ to $T-1$}
    \item Recede horizon $\mu^{t_{\text{mpc}}+1}_t = \mu^{t_{\text{mpc}}+1}_{t+1}$;
    \ENDFOR
    \ENDFOR
\end{algorithmic}
\label{alg: Adaptive_SS}
\end{algorithm}

\section{Algorithm}
\label{sec: algo}

In this section, we present the adaptive feedback stochastic search algorithm for opinion dynamics control implemented in \ac{MPC} fashion, as shown in \cref{alg: Adaptive_SS}.

At initial time, the feedback policy is initialized as $\calN(0, \Sigma)$, and the actuation probabilities are initialized as $p^0=\zeta$. At each MPC step, the current state is used to initialize the $M$ trajectory rollouts. For each rollout, the actuation variable $a^{i,m}$ is sampled from a Bernoulli distribution with $a^{i,m}=1$ indicating that the agent is actuated. The feedback policy parameter $K_t$ and noise $w_t^i$ are then sampled for each timestep and used to propagate the dynamics over the planning horizon. A cost is computed for each state-control trajectory rollout and used to update the feedback policy parameter and actuation probability distributions. Note that the actuation probability update \eqref{eq: bernoulli_numerical_update} can lead to changing active agent set size. To avoid this issue, a normalization step is performed after the update. This ensures that the expected active agent set size remains constant as $\Eb[\frac{1}{N}\sum_{i=1}^N p_i^k]=\xi$ throughout optimization. The active agent set for the current MPC step is decided based on a new set of actuation indicators sampled from the updated actuation probabilities. The first updated feedback policy is then applied to the current state. Finally, we recede the optimization horizon and re-initialize the feedback policy at the last timestep.

Note that the outlined algorithm performs a single iteration of optimization update at each MPC step. In addition, only the first policy step is executed every time. The proposed algorithm can be straightforwardly generalized to multiple optimization iterations and policy execution. On the other hand, a single realization of the dynamics noise is sampled for each rollout. This corresponds to the stochastic approximation scheme \cite{robbins1951stochastic} from stochastic optimization literature. A batched version of the algorithm can be derived similar to \cite{wang2021adaptive}.

\begin{figure}[t]
    \centering
    \includegraphics[width=0.7\linewidth]{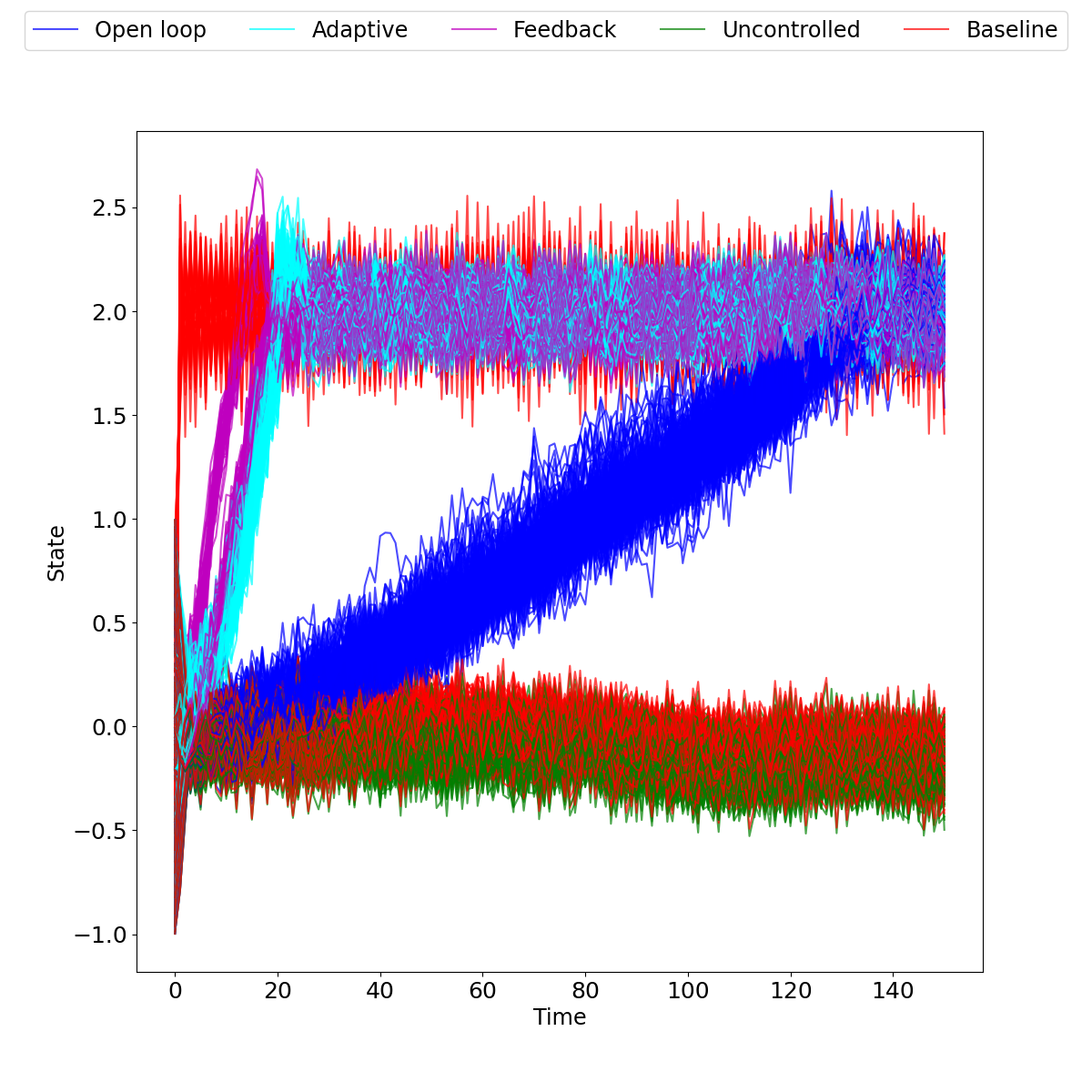}
    \caption{State trajectory comparison of different policies for the opinion steering setup with 25\% actuation percentage.}
    \label{fig: steering_comparison}
\end{figure}

\section{Simulation}
\label{sec: simulation}

\begin{figure*}[t!]
    \centering
    \begin{subfigure}{0.325\textwidth}
    \includegraphics[width=\textwidth]{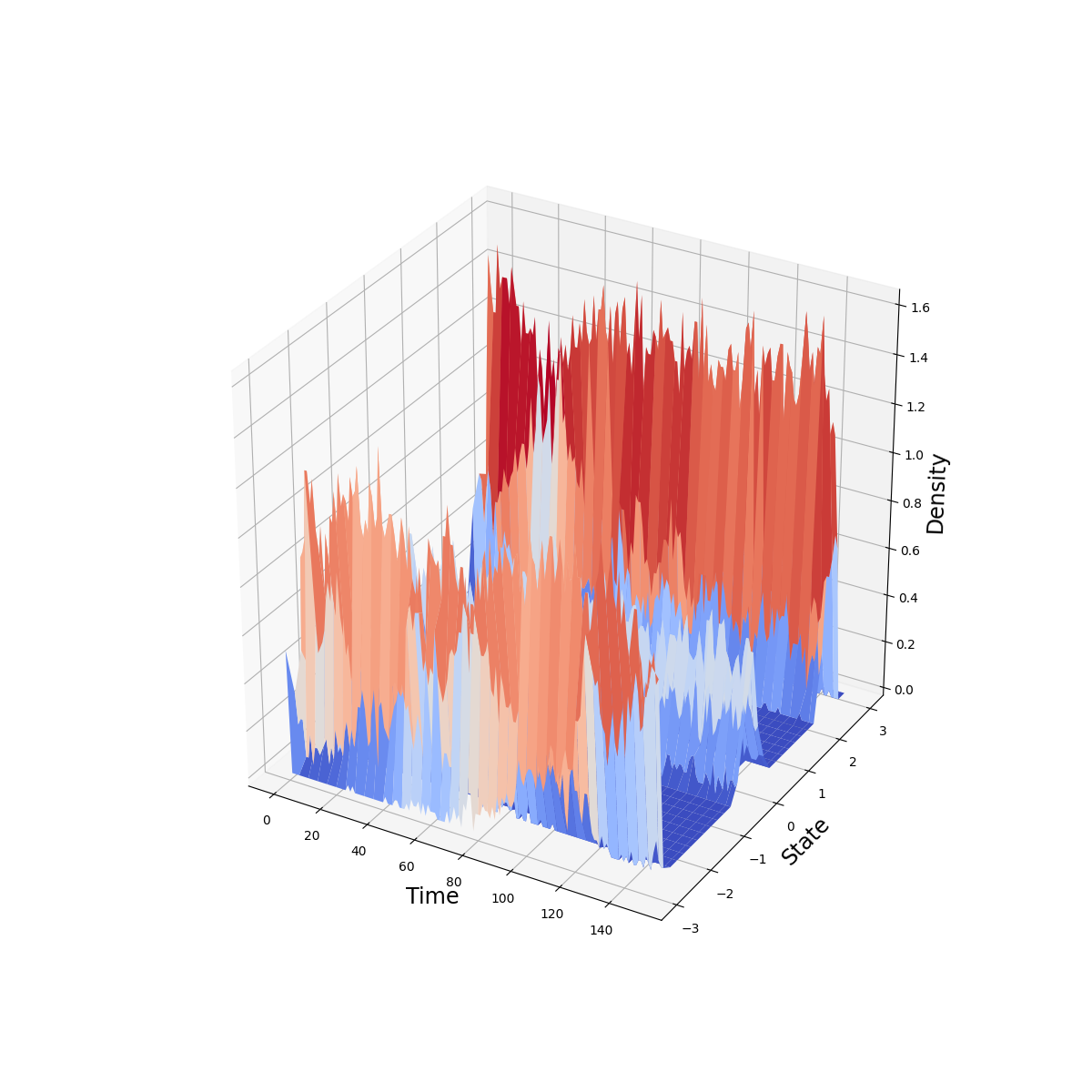}
    \caption{Baseline}
    \label{fig: manual}
\end{subfigure}
\begin{subfigure}{0.325\textwidth}
    \includegraphics[width=\textwidth]{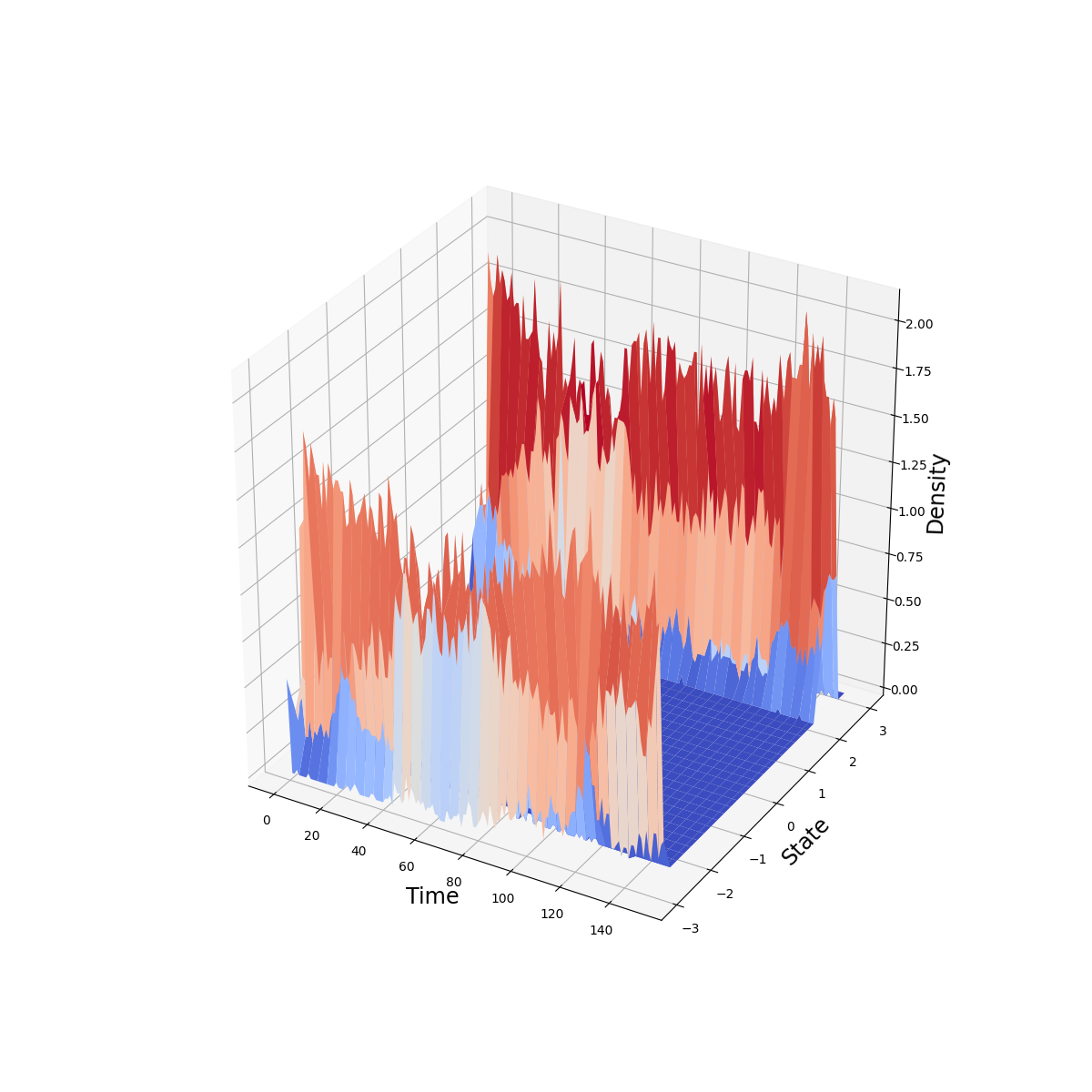}
    \caption{Uncontrolled}
    \label{fig: uncontrolled}
\end{subfigure}
\begin{subfigure}{0.325\textwidth}
    \includegraphics[width=\textwidth]{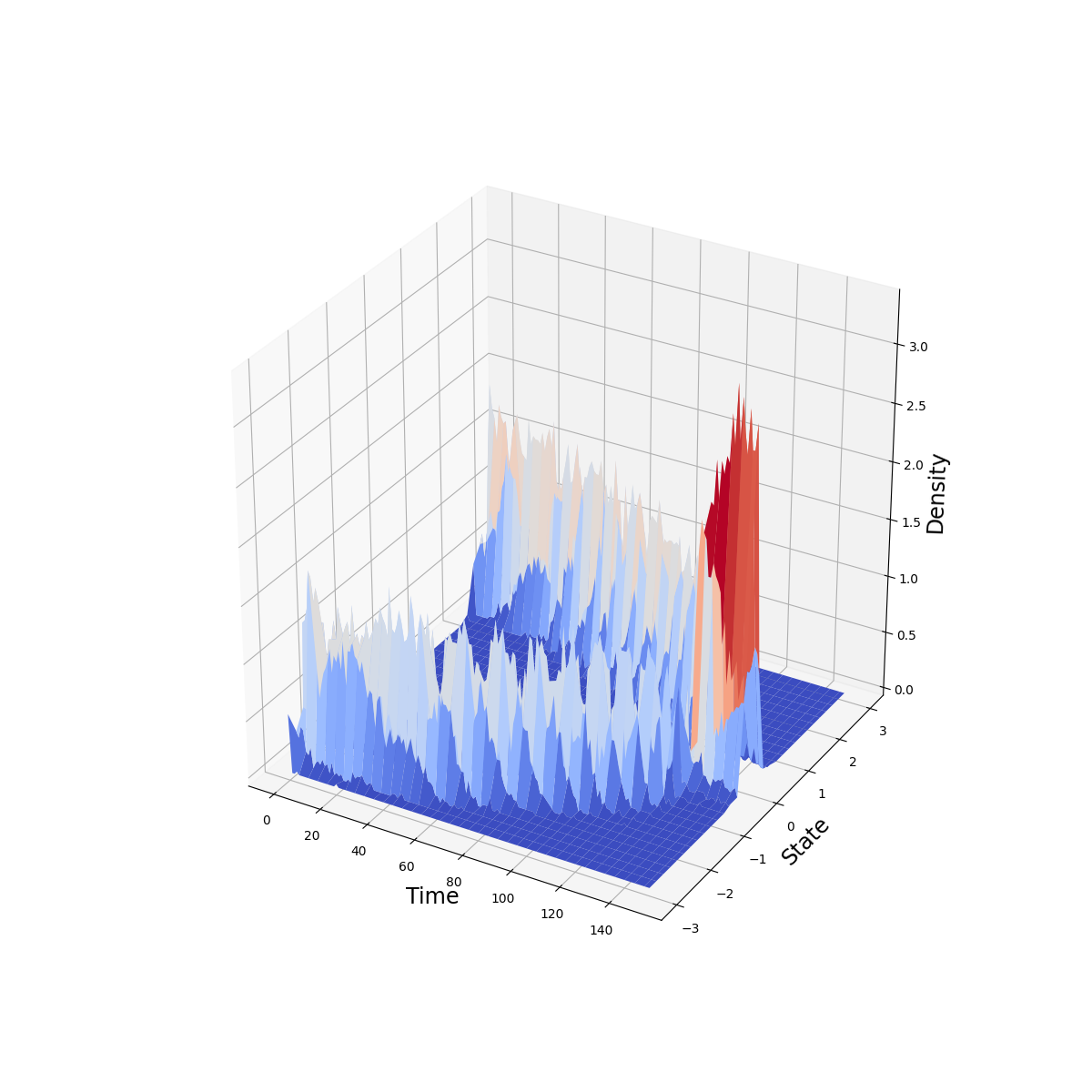}
    \caption{Open loop policy}
    \label{fig: openloop}
\end{subfigure}
\begin{subfigure}{0.325\textwidth}
    \includegraphics[width=\textwidth]{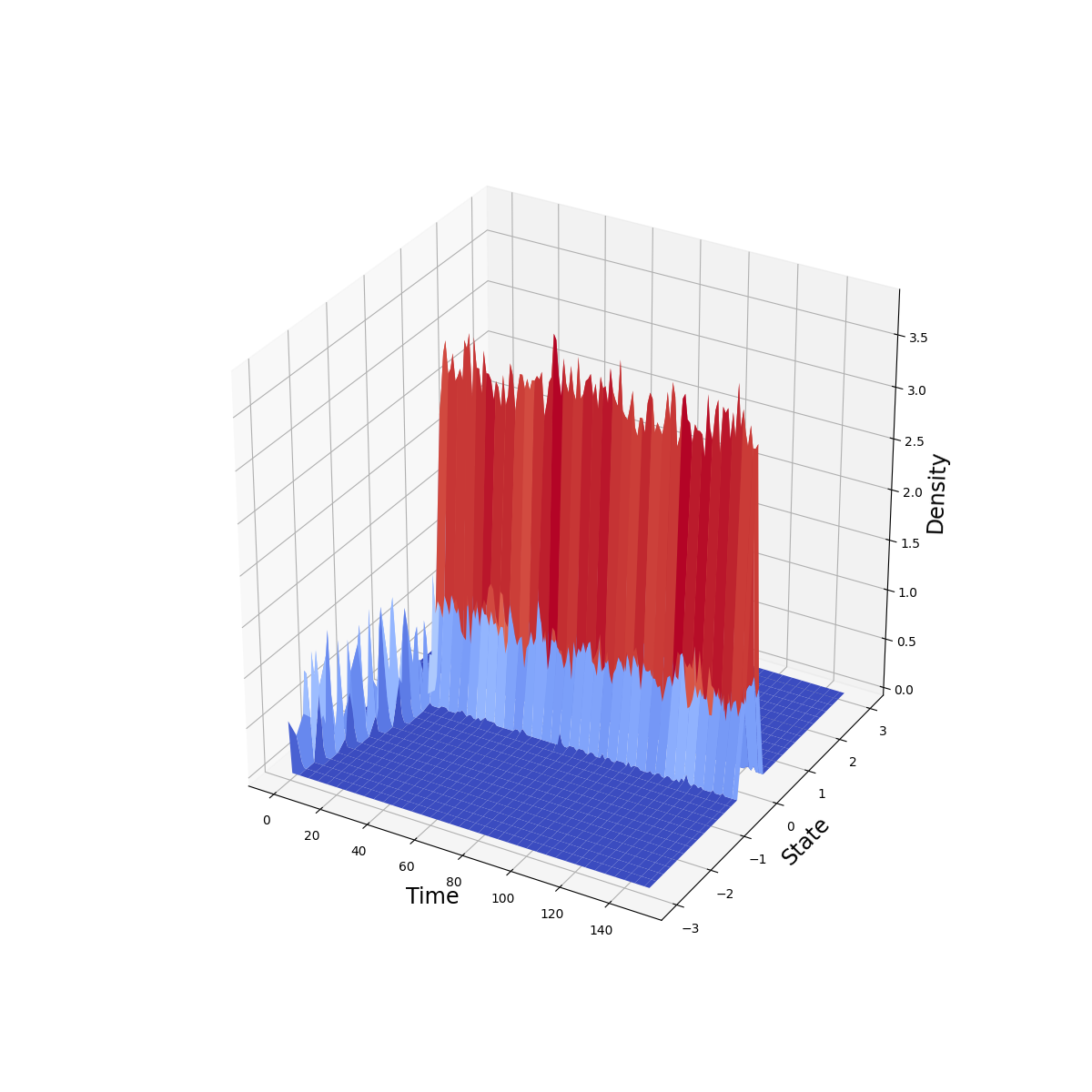}
    \caption{Feedback policy}
    \label{fig: feedback}
\end{subfigure}
\begin{subfigure}{0.325\textwidth}
    \includegraphics[width=\textwidth]{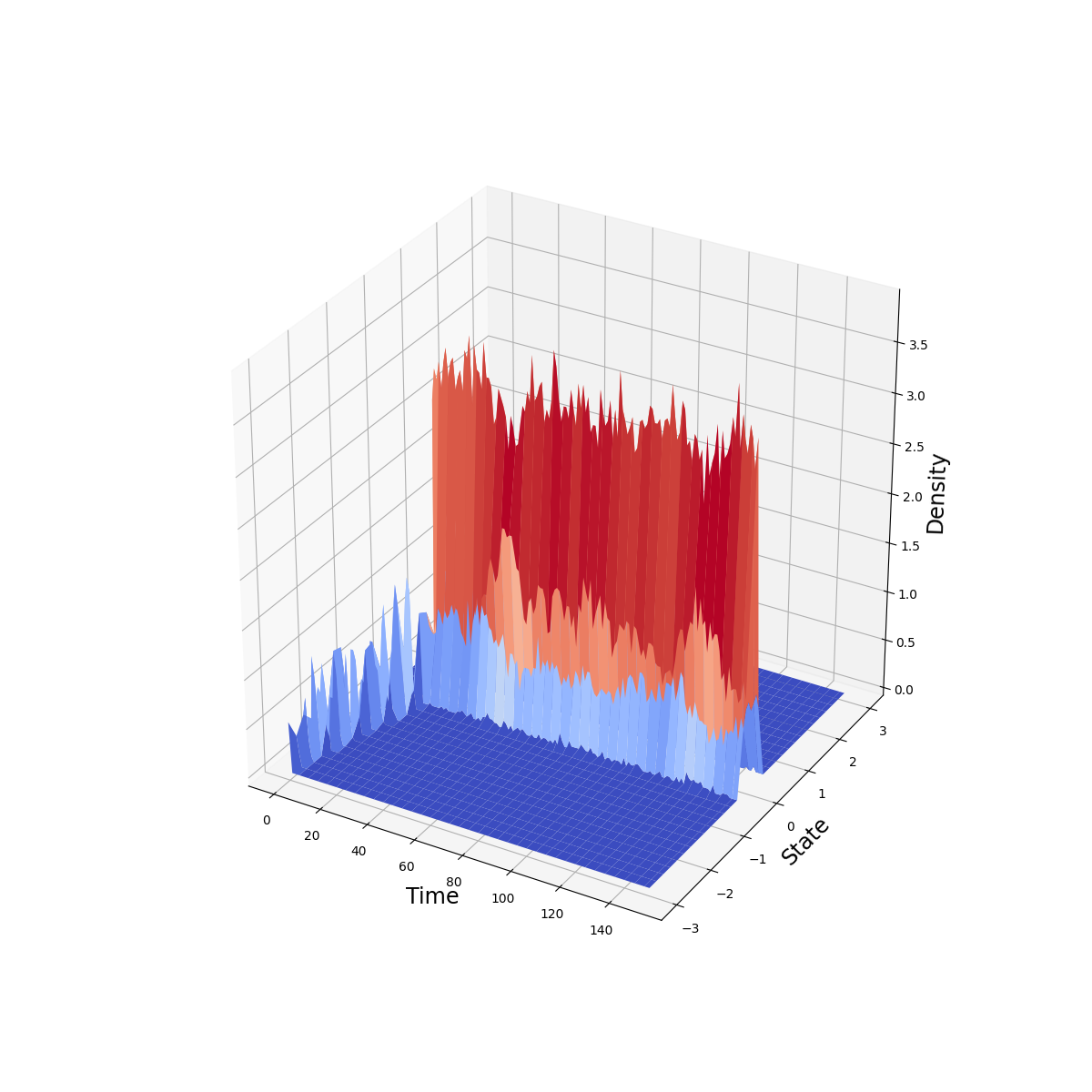}
    \caption{Adaptive feedback policy}
    \label{fig: adaptive}
\end{subfigure}
\begin{subfigure}{0.325\textwidth}
    \includegraphics[width=\textwidth]{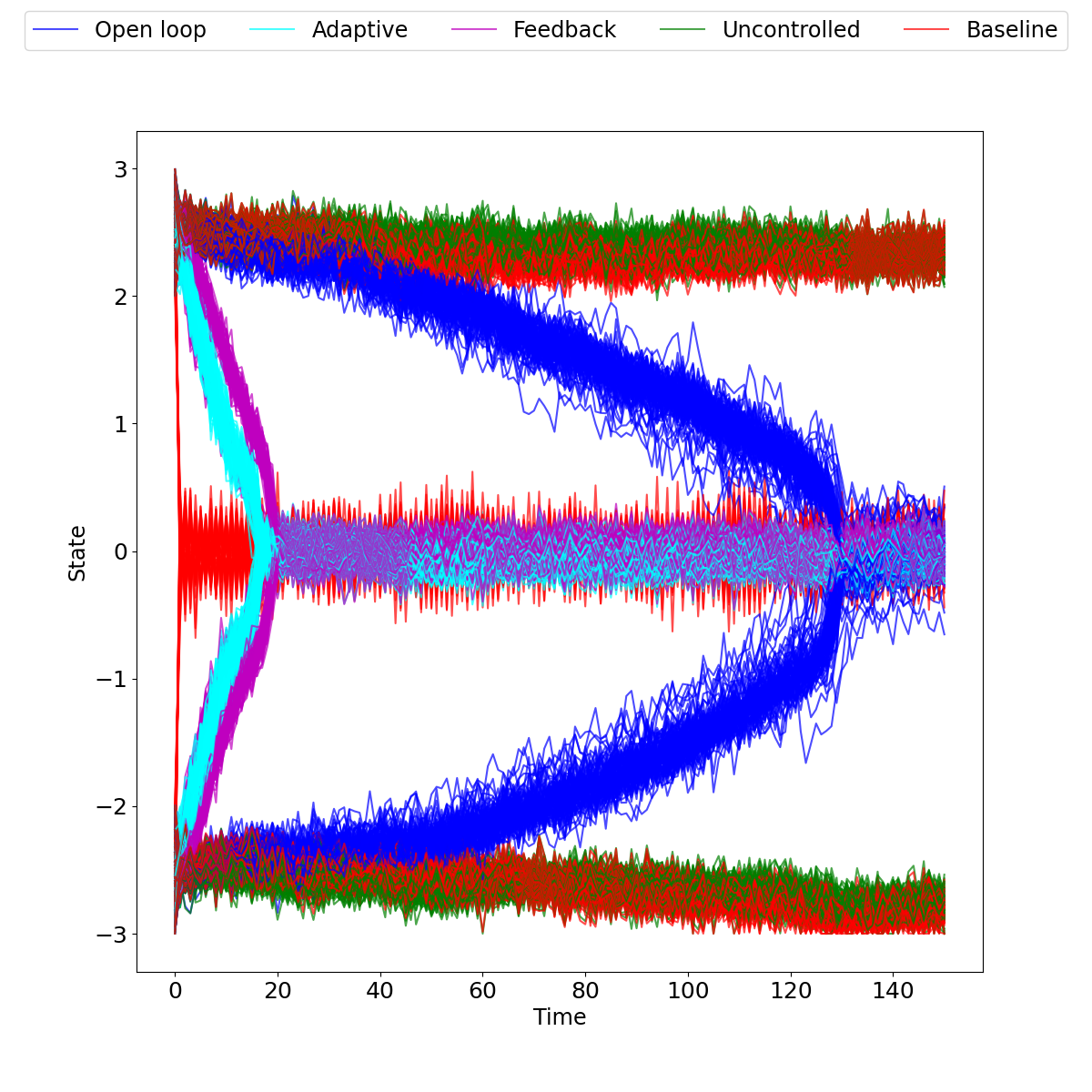}
    \caption{Comparison}
    \label{fig: comparison}
\end{subfigure}
    \caption{Comparison of different policies for the unifying polarized opinion setup. (a) - (e) are the population density plot over time, and (e) compares the state trajectory of each agent.}
    \label{fig: concensu_simulation}
\end{figure*}

In this section we present the simulation results of the proposed algorithm for opinion dynamics control in two different scenarios. We compare the performance of the following policies:
\begin{itemize}
\item \textbf{Open loop:} $u^i_t = \phi^i_t, \forall i\in I_A$
\item \textbf{Feedback:} $u^i_t = K_t x^i_t + k_t, \phi_t = \{K_t, k_t\}, \forall i\in I_A$
\item \textbf{Adaptive feedback:} $u^i_t = K_t x^i_t, \phi_t = K_t, \forall a_i=1$
\end{itemize}
against a hand-designed feedback control of $u^i_t=\zeta-x^i_t$ as the baseline. The baseline control has the flavor of pinning control as a function of the distance to the desired state, but it is unaware of the neighborhood. In addition, the control-free dynamics is included as a reference for the equilibrium population distribution.

\subsection{Simulation Setup}
In both scenarios, we set the population size $N$ as 200, active portion $\xi$ as $0.25$, susceptibility to local influence $\alpha$ as $0.8$, noise level $\sigma$ as $0.1$, and neighborhood size $\epsilon$ as $1$. The state of each agent is bounded in $[-3,3]$. The active agents are randomly generated based on the active agent set size. For the cost function, state cost coefficient $q=5$ and control cost coefficient $r=0.1$ are used. For the stochastic search controller, the planning horizon $T$ is $10$ and the algorithm is receded over $T_{\text{mpc}}=150$ MPC steps. During optimization, 500 samples are generated and the sigmoid shape function is used to transform the cost function with $\psi, \lambda, y_{lb}$ values that correspond to a $0.1$ elite threshold. A constant step size $\beta$ of $1$ is used.

\subsection{Opinion Dynamics Steering}
We first demonstrate the capability of the proposed framework on a unimodal steering task. The initial state of each agent is sampled from a uniform distribution in $[-1, 1]$. The goal is to steer all agents to $+2$. From the state trajectories comparison plot in \cref{fig: steering_comparison}, we can observe that the uncontrolled population stays around the neutral position. The baseline policy steers the active agents to the target position, but the passive agents remain around neutral position. While all three stochastic search policies manage to steer the population to the target location, the open loop policy reaches the desired state much more slowly than the feedback and adaptive feedback policies. Note that the feedback policy trajectories are split, with the passive agents trailing the active one, while the adaptive policy avoids the lag by adjusting the active agent set.

\subsection{Unifying Polarized Opinion}
We also test the different policies in a unifying polarized opinion setup. The initial state distribution in this case consists of two modes: $\text{Uniform}(2,3)$ and $\text{Uniform}(-3,-2)$. The goal is to bring both modes to the neutral position. Since the neighborhood size $\epsilon=1$ and noise level $\sigma=0.1$, the agents in each mode are only influenced by other agents in the same mode. Therefore, the uncontrolled population stays in the polarized bimodal distribution, as shown in \cref{fig: uncontrolled}. On the other hand, similar to the previous simulation setup, the baseline policy steers the active agents to the neutral position right away, but they lose the influence on the passive agents from the two modes, as shown in \cref{fig: manual}. From \cref{fig: openloop} to \cref{fig: adaptive}, we observe that all three policies manage to bring all agents to the neutral position, with open loop policy being the slowest and adaptive feedback policy being the fastest. A comparison of all policies against the baseline and uncontrolled case is included in \cref{fig: comparison}.

\begin{figure}[t]
    \centering
    \includegraphics[width=0.75\linewidth]{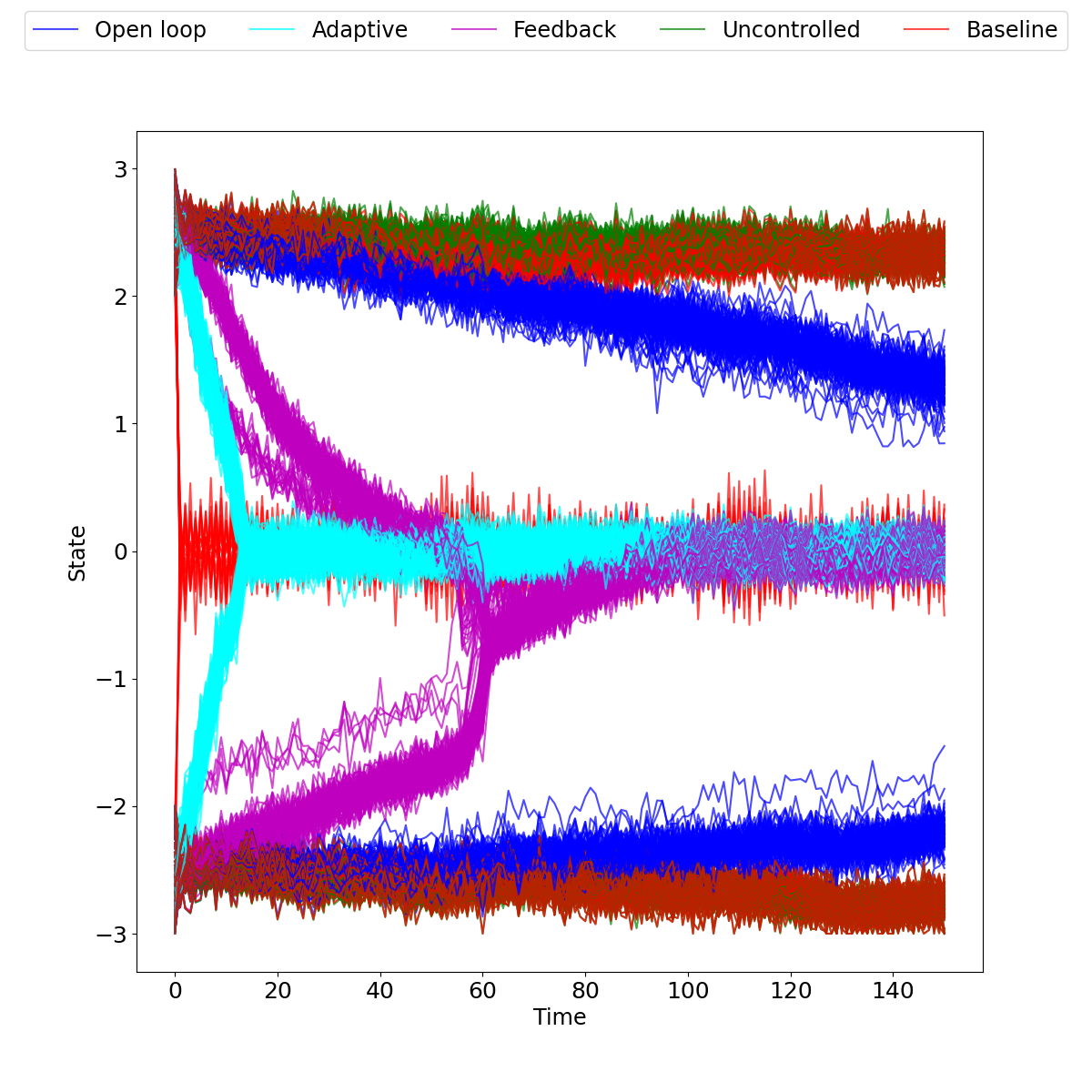}
    \caption{State trajectory comparison of different policies for the unifying polarized opinion setup with 10\% actuation percentage.}
    \label{fig: lowcontrol_comparison}
\end{figure}

An additional comparison is performed for the same setup but with active portion $\xi$ and $0.1$. From \cref{fig: lowcontrol_comparison}, we observe that the adaptive feedback policy performs similarly to \cref{fig: comparison} despite working with a much smaller active agent set. On the other hand, it takes the open loop and feedback policy much longer to steer the distribution. Note that the feedback policy trajectories indicate that the number of active agent in the negative mode is smaller than that of the positive mode, leading to worse steering performance for the mode.

\subsection{Discussions}

While open loop policy has been widely popular for sampling based MPC \cite{williams2017information, wang2021adaptive, wang2017variational}, the results of this paper showcase the power of feedback policy. Despite the implicit feedback through the MPC setup, explicit policy feedback parameterization results in faster convergence in both scenarios. In addition, the adaptive feedback policy can handle very small active portion with comparable performance.

\section{Conclusions and Future Work}
\label{sec: conclusion}

In this work, we presented the stochastic search framework for opinion dynamics steering where a subset of all agents are actuated. We derived the optimization update for an open loop, feedback, and a novel adaptive feedback policy that adjusts the active agent set. All three policies were tested on 2 opinion steering scenarios with 200 agents against a hand designed baseline policy and the uncontrolled dynamics. The results showcased the capability of the framework and the robustness of the adaptive feedback policy with respect to the active agent set size. As future work, we would like to scale the algorithm to large scale experiments through distributed optimization techniques such as ADMM. In addition, we would like to extend the framework to more complex opinion dynamics/social network models and design novel policies for the model.

%%%%%%%%%%%%%%%%%%%%%%%%%%%%%%%%%%%%%%%%%%%%%%%%%%%%%%%%%%%%%%%%%%%%%%%%%%%%%%%%
\section*{ACKNOWLEDGMENTS}

This work is supported by the DoD Basic Research Office Award HQ00342110002.

%%%%%%%%%%%%%%%%%%%%%%%%%%%%%%%%%%%%%%%%%%%%%%%%%%%%%%%%%%%%%%%%%%%%%%%%%%%%%%%%

\bibliographystyle{unsrt}
\bibliography{reference}

% \begin{thebibliography}{99}

% \bibitem{c1}
% J.G.F. Francis, The QR Transformation I, {\it Comput. J.}, vol. 4, 1961, pp 265-271.

% \bibitem{c2}
% H. Kwakernaak and R. Sivan, {\it Modern Signals and Systems}, Prentice Hall, Englewood Cliffs, NJ; 1991.

% \bibitem{c3}
% D. Boley and R. Maier, "A Parallel QR Algorithm for the Non-Symmetric Eigenvalue Algorithm", {\it in Third SIAM Conference on Applied Linear Algebra}, Madison, WI, 1988, pp. A20.

% \end{thebibliography}

\end{document}

%% file: packages.tex
\usepackage{hyperref}
\hypersetup{colorlinks = true,
            linkcolor = blue,
            urlcolor  = blue,
            citecolor = blue,
            anchorcolor = blue}

\usepackage{url}

\usepackage[utf8]{inputenc} % allow utf-8 input
\usepackage[T1]{fontenc}    % use 8-bit T1 fonts
\usepackage{url}            % simple URL typesetting
\usepackage{booktabs}       % professional-quality tables
\usepackage{amsfonts}       % blackboard math symbols
\usepackage{nicefrac}       % compact symbols for 1/2, etc.
\usepackage{microtype}      % microtypography
\usepackage{algorithmic}
\usepackage{algorithm}
\usepackage{amsmath}
\usepackage{amssymb}
\usepackage{bm}
\usepackage{mathtools}
\usepackage{mathrsfs}
\usepackage{subcaption}
\usepackage{graphicx} % for pdf, bitmapped graphics files
\usepackage{caption}
\usepackage[nolist]{acronym}
\usepackage{bbm}
\usepackage{xcolor}
\usepackage{wrapfig}
\usepackage{siunitx}
\usepackage{multirow}
\usepackage{todonotes}

\usepackage[symbol]{footmisc}
\usepackage{float}
\usepackage{cleveref}
\usepackage{caption}

\usepackage{tikz}

%% file: commands.tex
\newcommand{\T}{^\top}

\newcommand{\calX}{{\cal X}}

\newcommand{\calN}{{\cal N}}

\newcommand{\argmax}{\operatornamewithlimits{argmax}}
\newcommand{\argmin}{\operatornamewithlimits{argmin}}

\newcommand{\Rb}{\mathbb{R}}

\newcommand{\Eb}{\mathbb{E}}

\begin{acronym}
\acro{CEM}{Cross Entropy Method}
\acro{GASS}{Gradient-based Adaptive Stochastic Search}
\acro{MPPI}{Model Predictive Path Integral}
\acro{SA}{Stochastic Approximation}
\acro{DNN}{Deep Neural Network}
\acro{DP}{Dynamic Programming}
\acro{FBSDE}{Forward-Backward Stochastic Differential Equation}
\acro{FSDE}{Forward Stochastic Differential Equation}
\acro{BSDE}{Backward Stochastic Differential Equation}
\acro{LSTM}{Long-Short Term Memory}
\acro{FC}{Fully Connected}
\acro{DDP}{Differential Dynamic Programming}
\acro{HJB}{Hamilton-Jacobi-Bellman}
\acro{PDE}{Partial Differential Equation}
\acro{PI}{Path Integral}
\acro{NN}{Neural Network}
\acro{GPs}{Gaussian Processes}
\acro{SOC}{Stochastic Optimal Control}
\acro{RL}{Reinforcement Learning}
\acro{MPOC}{Model Predictive Optimal Control}
\acro{IL}{Imitation Learning}
\acro{RNN}{Recurrent Neural Network}
\acro{FNN}{Feed-forward Neural Network}
\acro{Single FNN}{Single Feed-forward Neural Network}
\acro{DL}{Deep Learning}
\acro{SGD}{Stochastic Gradient Descent}
\acro{SDE}{Stochastic Differential Equation}
\acro{2BSDE}{second-order BSDE}
\acro{2FBSDE}{second-order FBSDE}
\acro{OCE}{Optimized Certainty Equivalents}
\acro{CVaR}{Conditional Value-at-Risk}
\acro{VaR}{Value-at-Risk}
\acro{DRO}{Distributionally Robust Optimization}
\acro{PMP}{Pontryagin Maximum Principle}
\acro{VI}{Variational Inference}
\acro{ELBO}{Evidence Lower BOund}
\acro{EM}{Expectation Maximization}
\acro{SVGD}{Stein Variational Gradient Descent}
\acro{RKHS}{Reproducing Kernel Hilbert Space}
\acro{SME}{Stochastic Modified Equation}
\acro{BO}{Bayesian Optimization}
\acro{SS}{Stochastic Search}
\acro{ARA}{Absolute Risk Aversion}
\acro{QGASS}{Quantum \ac{GASS}}
\acro{MPC}{Model Predictive Control}
\acro{GPU}{Graphic Processing Unit}
\acro{VO}{Variational Optimization}
\acro{fGn}{fractional Gaussian noise}
\end{acronym}

\newtheorem{innercustomgeneric}{\customgenericname}
\providecommand{\customgenericname}{}
\newcommand{\newcustomtheorem}[2]{%
  \newenvironment{#1}[1]
  {%
   \renewcommand\customgenericname{#2}%
   \renewcommand\theinnercustomgeneric{##1}%
   \innercustomgeneric
  }
  {\endinnercustomgeneric}
}

\newcustomtheorem{customthm}{Theorem}
\newcustomtheorem{customlemma}{Lemma}